\documentclass{amsart}
\usepackage{palatino}
\usepackage{amsthm, amsmath}
\usepackage{amssymb} 
\usepackage{amscd}
\usepackage[margin=3cm]{geometry}
\usepackage[breaklinks,bookmarksopen,bookmarksnumbered]{hyperref}
\usepackage[all]{xy}

\theoremstyle{plain}

\newtheorem{prop}{Proposition}

\theoremstyle{remark}
\newtheorem{rem}{Remark}
\newtheorem*{ex}{Example}

\newcommand\pr{\noindent\textit{Proof} : }

\newcommand\rond{\kern 1pt{\scriptstyle\circ}\kern 1pt}

\newcommand\Ker{\operatorname{Ker}}

\newcommand\Pic{\operatorname{Pic}}

\newcommand\Tr{\operatorname{Tr}}
\newcommand\fx{\varphi^{}_X}
\newcommand\ft{\varphi^{}_S}

\newcommand\Z{\mathbb{Z}}

\renewcommand\P{\mathbb{P}}

\renewcommand\O{\mathcal{O}}

\def\qfl#1{\buildrel {#1}\over {\longrightarrow}}

\newcommand\iso{\vbox{\hbox to .8cm{\hfill{$\scriptstyle\sim$}\hfill}
\nointerlineskip\hbox to .8cm{{\hfill$\longrightarrow $\hfill}} }}

\begin{document}
\title[A tale of two surfaces]{A tale of two surfaces}
\author[Arnaud Beauville]{Arnaud Beauville}
\address{Laboratoire J.-A. Dieudonn\'e\\
UMR 7351 du CNRS\\
Universit\'e de Nice\\
Parc Valrose\\
F-06108 Nice cedex 2, France}
\email{arnaud.beauville@unice.fr}
 
\begin{abstract}
We point out a link between two  surfaces which have appeared recently in the literature: the surface of cuboids and the Schoen surface. Both  give rise to a surface  with $q=4$, whose  canonical map  is 2-to-1 onto a complete intersection of 4 quadrics  in $\P^6$ with 48 nodes.
\end{abstract}

\maketitle 
\section{Introduction}
The aim of this note is to point out a link between two  surfaces which have appeared recently in the literature: the \emph{surface of cuboids} (\cite{ST}, \cite{vL}) and the surface (actually a family of surfaces) discovered by Schoen \cite{S}. We will show that both surfaces give rise to a surface $X$  with $q=4$, whose  canonical map  is 2-to-1 onto a complete intersection of 4 quadrics  $\Sigma \subset \P^6$ with 48 nodes. In the first case (\S 2)  $X $ is a quotient $ (C\times C')/ (\Z/2)^2  $, where $ C $ and $ C' $ are genus $ 5 $ curves with a free action of $ (\Z/2)^2 $.
In the second case (\S 3), $X$ is a double \'etale cover of the Schoen surface.

When the canonical map of a surface $X$ of general type has degree $>1$ onto a surface,  that surface  either  has $p_g=0$ or  is itself canonically embedded (\cite{B}, Th. 3.1). Our surfaces $X$ provide one more example of the latter case, which is rather exceptional (see \cite{CPT} for a list of the examples known so far).

\medskip	
\section{The surface of cuboids and its deformations}
\medskip	

In $ \P^4 $, with coordinates $ (x,y;u,v,w)  $, we consider the curve $ C $ given by
\begin{equation}
 u^2=a(x,y) \quad,\quad v^2=b(x,y) \quad,\quad w^2=c(x,y)   
\end{equation}
where $ a,b,c $ are quadratic forms in $ x,y $. 
We assume that the zeros $ \{p'_a, p''_a; p'_b, p''_b;p'_c, p''_c\}$ of $ a,b,c $ form a set $ Z\subset \P^1 $ of 6 distinct points.
Then $ C $ is a smooth curve of genus $ 5 $, canonically embedded. It is preserved by 
the group  $ \Gamma_+ \cong (\Z/2)^3  $ which acts on $ \P^4 $ by changing the signs of $ u,v,w $. Let $\Gamma\subset\Gamma_+$ be the subgroup  (isomorphic to $(\Z/2)^2  $) which changes an even number of signs. It acts freely on $ C $, so the quotient curve $ B:=C/\Gamma $ has genus 2.  The subring of $ \Gamma $-invariant elements in 
$ \oplus\ H^0(C,K_C^{n})  $ is generated by $ x,y $ and $ z:=uvw $, with the relation  $ z^2=abc $; thus $ B $ is 
the double cover of $ \P^1 $ branched along $ Z $. 

Let  $ JB_2$ be the group of 2-torsion line bundles on $ B $ (isomorphic to $  (\Z/2)^4 $). The $ \Gamma $-covering $\pi: C\rightarrow B $ corresponds to a subgroup 
$ \cong (\Z/2)^2$ of $ JB_2 $, namely the kernel of $ \pi^*:JB\rightarrow JC $. Since the divisor 
 $ \pi^*(p'_a+p''_a)  $ is cut out on $ C $ by the canonical divisor $ u=0 $, we have $ \pi^* (p'_a-p''_a)\sim 0$, and similarly for $ b $ and $ c $; thus 
$\Ker \pi^*= \{0,p'_a-p''_a,p'_b-p''_b,p'_c-p''_c\} $. This is a \emph{Lagrangian} subgroup of $ JB_2 $ for the Weil pairing \cite{M2}; conversely, any Lagrangian subgroup of $ JB_2 $ is of that form. Thus the curves $ C $ we are considering are exactly the $  (\Z/2)^2 $-\'etale covers of a curve $ B $ of genus 2 associated to a Lagrangian subgroup of $ JB_2 $. In particular they form a 3-dimensional family.

The group  $ \Gamma_+/\Gamma\cong \Z/2 $ acts on $ B=C/\Gamma $ through the hyperelliptic involution, so $ C/\Gamma_+ $ is isomorphic to $ \P^1 $.

\begin{prop}
Let $ C,C' $ be two genus $ 5 $ curves of type $ (1)  $, and let $ X $ be the quotient of  $ C\times C' $ by  the diagonal action of $ {\Gamma}\cong (\Z/2)^2$. 

$ 1) $ $ X $ is a minimal surface of general type with $ q=4 $, $ p_g=7 $, $ K^2=32 $.

$ 2) $ The involution $ i_X $ of $ X $  defined by the action of $ \Gamma_+ /\Gamma \cong \Z/2$ has $ 48 $ fixed points. The canonical map $\fx :X\rightarrow \P^6$ factors through $i_X$, and induces an isomorphism of $X/i_X$ onto a complete intersection of $4$ quadrics in $\P^6$ with $48$ nodes. 
\end{prop}
\pr The computation of the numerical invariants of $ X $ is straightforward.

Let us denote by $ (x',y';u',v',w')  $ the coordinates on $ C' $, and by $ a',b',c' $ the corresponding quadratic forms. A basis of the space $ H^0(X,K_X)= \bigl( H^0( C,K_C)\otimes H^0(C',K_{C'})  \bigr)^{\Gamma}$ is given by the elements
\[ X= x \otimes x'\quad Y= x\otimes y' \quad Z= y \otimes x'\quad T=y \otimes y'\quad U=u \otimes u'\quad V=v \otimes v'\quad  W= w \otimes w'\ . \]
They satisfy the relations
\[ XT-YZ=0 \quad,\quad U^2 = A(X,Y,Z,T) \quad,\quad  V^2 = B(X,Y,Z,T) \quad,\quad W^2 = C(X,Y,Z,T) \ ,\]
where $ A,B,C $ are quadratic forms satisfying $ A(X,Y,Z,T)=a(x,y)\otimes a(x',y')   $ and the analogous relations for $ B $ and $ C $.

Let $ \Sigma $ be the surface defined by these 4 quadratic forms, and let $ \varphi:X\rightarrow \Sigma $ be the induced map.  We have 
$ \varphi\rond i_X=\varphi $, so $ \varphi $ induces a map $ \bar{\varphi}$ from $X/i_X= (C\times C' )/\Gamma_+$ into $ \Sigma $. We consider the commutative diagram
\[ \xymatrix{(C\times C' )/\Gamma_+ \ar[rr]^{\bar{\varphi}}\ar[dr]_p&& \Sigma\ar[dl]^q\\
&Q\cong \P^1\times\P^1&
} \]

\noindent where $ p : (C\times C' )/\Gamma_+\rightarrow (C/\Gamma_+)\times (C' /\Gamma_+) $  is the quotient map by $ \Gamma_+ $, and $ q $ the projection \allowbreak$ (X,Y,Z,T;U,V,W)\allowbreak\mapsto  (X,Y,Z,T)  $. The group $(\Z/2)^3 $ acts on 
$ \Sigma $ by changing the signs of $ (U,V,W)  $; then  $ \bar{\varphi} $ is an equivariant map of $(\Z/2)^3 $-coverings, hence an isomorphism.

\smallskip
It remains to show that $ i_X $ has 48 fixed points. These fixed points are the images $ (\mathrm{mod.}\ \Gamma)  $ of the points of $ C\times C' $ fixed by one of the elements  of $ \Gamma_+\smallsetminus \Gamma $. Such an element changes the sign of one of the coordinates $\ell= u,v $ or $ w $,  hence fixes the 64 points $ (m,m')  $ of $ C\times C' $ with $\ell(m)=\ell(m')  =0  $. This gives $ (3\times 64)/4=48 $ fixed  points in $ X $.\qed 

\medskip
\begin{ex}
Let us take for $ C$ and $ C' $ the curve
\[ u^2=xy \quad,\quad v^2=x^2-y^2\quad,\quad w^2=x^2+y^2\ . \]
The set of zeros of $ a,b,c $ is $ \{0,\infty,\pm 1,\pm i\} $, so the genus 2 curve $ B $ is given by $ z^2=x(x^4-1)  $.

We get for $ \Sigma $ the following equations  :
\[ XT=YZ= U^2 \quad , \quad  V^2 = X^2-Y^2-Z^2+T^2 \quad,\quad W^2 = X^2+Y^2+Z^2+T^2\ ;  \]
or, after the linear change of variables $ X=\mathsf{x+t} $, $ T=\mathsf{t-x} $, $ Y= \mathsf{y+iz} $, $ Z=\mathsf{y-iz}$, $ U=\mathsf{u} $, $ V= 2\mathsf{v}$, $ W=2\mathsf{w} $:
\[ \mathsf{t^2= x^2+y^2+z^2 \quad , \quad u^2= y^2+z^2 \quad , \quad v^2=x^2+z^2 \qquad w^2=   x^2+y^2} \ .  \] 
These are the equations of the \emph{surface of cuboids}, studied in  \cite{ST}, \cite{vL}. It encodes the relations in a cuboid (= rectangular box) between the sides $ \mathsf{x,y,z}$, the diagonals of the faces $ \mathsf{u,v,w}$, and the big diagonal $ \mathsf{t}$. Thus the surface of cuboids belongs to a 6-dimensional family of intersection of 4 quadrics in $ \P^6 $ with 48 nodes.
\end{ex}
\begin{rem}
The surfaces $ X $ fit into a tower of $ (\Z/2)^2  $-\'etale coverings :
\[ C\times C' \longrightarrow X \qfl{r} B\times B' \ .\]
The abelian covering  $ r $ is the pull back of a $ (\Z/2)^2  $-\'etale covering of $ JB\times JB' $ :
\[ \xymatrix@M=7pt{X\ar@{^{(}->}[r]^\alpha\ar[d]^r & A\ar[d] \\
B\times B' \ar@{^{(}->}[r] & JB\times JB'\ .
}  \]

The abelian variety $ A $  is the Albanese variety of $ X $, and $ \alpha $ is the Albanese map. Since the quotient $ X/i_X $ is regular, $ i_X $ acts as $ (-1)  $ on the space $ H^0(X,\Omega^1_X)  $; therefore 
if we choose $ \alpha $ so that it maps a fixed point of $ i_X $ to $ 0 $, $ i_X $ is induced by $ (-1_A)  $.
\end{rem}

\medskip	
\section{The Schoen surface}

The Schoen surfaces $S$ have been defined in \cite{S}, and studied in \cite{CMR}. A Schoen surface $ S $
is contained in its Albanese variety $A$; it has the following properties:

 a) $\ K_S^2=16\ ,\  p_g=5\ ,\ q=4\ $ (hence $ \chi (\O_S)=2 $);

b)  The canonical map $\ft :S\rightarrow \P^4$ factors through an involution  $i_S$  with 40 fixed points, and induces an isomorphism of $S/i_S$ onto the complete intersection of a quadric and a quartic in $\P^4$ with 40 nodes \cite{CMR}.

\smallskip
Since $ S/i_S $ is a regular surface,  $ i_S $ acts  as $ (-1)  $ on the space $ H^0(S,\Omega^1_S)  $. Therefore
if we choose the Albanese embedding $ S\hookrightarrow A $ so that it maps a fixed point of $ i_S $ to $ 0 $, $ i_S $ \emph{is induced by the involution} $ (-1_A)  $.
\medskip

Let $\ell $ be a line bundle of order 2 on $A$; we denote by $\pi :B \rightarrow A$ the corresponding \'etale double cover, and put $X:=\pi ^{-1}(S)$.  The restriction of $\ell$ to $S$, which we will still denote by $\ell$, is nontrivial (because the restriction map $\Pic^{\mathrm{o}}(A)\rightarrow \Pic^{\mathrm{o}}(S)$ is an isomorphism), hence $X$
 is connected.

\begin{prop}\label{q}
$X$ is a minimal surface of general type with $q=4$, $p_g=7$, $K_X^2=32$.
\end{prop}
\pr The formulas $K_X^2=32$ and $\chi (\O_X)=4$ are immediate; we must prove $q(X)=4$, that is, $H^1(S,\ell)=0$. 

By construction \cite{S} the Schoen surfaces fit into a flat family over the unit disk $\Delta $:

\[\xymatrix@M=7pt{ \mathcal{S} \ar@{^{(}->}[r]\ar[dr]& \mathcal{A}\ar[d]\\
&\Delta 
}\]where: 

$\bullet \ \mathcal{A}/\Delta$ is a smooth family of abelian varieties; 

$\bullet$ at a point $z\neq 0$ of $\Delta $, $\mathcal{S}_z$ is a Schoen surface, and $\mathcal{S}_z\hookrightarrow  \mathcal{A}_z$ is the Albanese embedding;

$\bullet\ \mathcal{A}_0=JC\times JC$ for a genus 2 curve $C$; $\mathcal{S}_0$ is the union of $JC$ embedded diagonally in 
$JC\times JC$, and of $C\times C\subset JC\times JC$ (we choose an Abel-Jacobi embedding $C\subset JC$). These two components intersect transversally along the diagonal $C\subset C\times C $. 

The  line bundle $\ell$ extends to a line bundle $\mathcal{L}$ of order 2 on $\mathcal{A}$. Let $\ell_0 $ be the restriction of $\mathcal{L}$ to $\mathcal{S}_0 $; we want to compute $H^1(\mathcal{S}_0,\ell_0 )$. 
We have an exact sequence
\begin{equation}\label{ex}
0\rightarrow\ell_0  \longrightarrow  \ell_{0\,|JC}\oplus \ell_{0\,|C\times C} \longrightarrow \ell_{0\,|C}\rightarrow 0\ .
\end{equation}
The line bundle $\mathcal{L}_0$ on $JC\times JC$ can be written $\alpha \boxtimes \beta $, where $\alpha $ and $\beta $ are 
2-torsion line bundles on $JC$, not both trivial; we use the same letters to denote  their restriction to $C$. The cohomology exact sequence associated to (\ref{ex}) gives

\begin{multline*}
H^0(JC, \alpha \otimes \beta )\oplus H^0(C\times C, \alpha \boxtimes \beta )\longrightarrow H^0(C, \alpha \otimes \beta ) \longrightarrow H^1(\mathcal{S}_0,\ell_0 )\qfl{u} \\ H^1(JC, \alpha \otimes \beta )\oplus H^1(C\times C, \alpha \boxtimes \beta ) \longrightarrow H^1(C, \alpha \otimes \beta )\ . 
\end{multline*}

The restriction map $H^0(JC,\alpha \otimes \beta)\rightarrow H^0(C,\alpha \otimes \beta)$ is surjective, so $u$ is injective. If $\alpha $ and $\beta $ are nontrivial, $H^1(C\times C, \alpha \boxtimes \beta )$ is zero, and the restriction map $H^1(JC, \alpha \otimes \beta )\rightarrow H^1(C, \alpha \otimes \beta )$ is injective, so $H^1(\mathcal{S}_0,\ell_0)=0$. If, say, $\beta $ is trivial, $H^1(JC, \alpha  )$ is zero and  the map $H^1(C\times C, \mathrm{pr}_1^*\alpha )\rightarrow H^1(C,\alpha )$ is bijective, hence  $H^1(\mathcal{S}_0,\ell_0)=0$ again.

By semi-continuity this implies $H^1(\mathcal{S}_z,\mathcal{L}_z)=0$ for $z$ general in $\Delta $, or equivalently $q(\widetilde{\mathcal{S}}_z)=q({\mathcal{S}}_z)=4$, where 
$\widetilde{\mathcal{S}}\rightarrow \mathcal{S} $ is the \'etale double covering defined by $\mathcal{L}$. But $q$ is a topological invariant, so this holds for all $z\neq 0$ in $\Delta $, hence $H^1(S,\ell)=0$.\qed

\medskip	
 The surface $X$ has a natural action of $(\Z/2)^2$, given by 
 the involution $i_X$ induced by $(-1_B)$ and the involution $\tau $ associated to the double covering $X\rightarrow S$, which is induced by a translation  of $B$. We want to determine how these involutions act on $H^0(X,K_X)$. The decomposition of $H^0(X,K_X)$ into eigenspaces for $\tau $ is
\[ H^0(X,K_X) \cong H^0(S,K_S)\oplus H^0(S, K_S\otimes  \ell)\ .\]
 By property b) above, $i_S$ acts trivially on $H^0(S,K_S)$.   It remains to study how it acts on $H^0(S, K_S\otimes \ell)$, or equivalently on $H^2(S,\ell)$. To define this action we  choose the  isomorphism $u : (-1_A)^*\ell\iso\ell$ over $A$ such that $u(0)=1$, and we consider the involutions $H^p(i_S,u) : H^p(S,\ell)\qfl{i_S^*} H^p(S,i_S^*\ell)\qfl{u_{|S}} H^p(S,\ell)$.

\begin{prop}\label{f}
There exist  line bundles
 $\ell$ of  order 2 on $A$  for which $i_S$ acts trivially on $H^2(S,\ell)$. In that case $i_X$ has $48$ fixed points.
\end{prop}

 \pr  We will denote by $A_2$ and $\hat{A}_2$ the  2-torsion subgroups of $A$ and its dual abelian variety $\hat{A}$, and similarly for $B$. The fixed point set of $i_S$ is $A_2\cap S$, and that of $i_X$ is $B_2\cap X$. 
 
 We apply the holomorphic Lefschetz formula to the automorphism $i _S$ of $S$ and the $ i_S $-linearization $u_{|S}:i_S^*\ell\rightarrow \ell$:
 \[  \sum_p (-1)^p \Tr H^p(i_S,u)  = \frac{1}{4} \sum_{a\in A_2\cap S} u(a) \ .
\]
 (At a point $a$ of $A_2$, $u(a):\ell_a\rightarrow \ell_a$ is the multiplication by a scalar, which we still denote $u(a)$.)

Let $a\in A_2$.  By \cite{M1}, property iv) p.\ 304, we have 
 $u(a)=(-1)^{\langle  a,\ell \rangle}$, where $\langle \ ,\  \rangle: A_2\times \hat{A}_2\rightarrow \Z/2$ is the canonical pairing. 
On the other hand, dualizing the exact sequence of $(\Z/2)$-vector spaces
\[ 0\rightarrow (\Z/2)\ell \rightarrow \hat{A}_2 \qfl{\hat \pi }  \hat{B}_2 \]
and using the canonical pairings we get an exact sequence
\[B_2\qfl{\pi } A_2 \qfl{\langle \ ,\ell \rangle}\Z/2\rightarrow 0\ .\]
Thus $u(a)=1$ if  $a\in \pi (B_2)$, and $u(a)=-1$ otherwise. For $ i=0 $ or $ 1 $, let $f_i$  be the number of points $a\in A_2\cap S$ with $ {\langle  a,\ell \rangle}=i $.
The right hand side of the Lefschetz formula is 
$\frac{1}{4}(f_0-f_1) $. 

\smallskip	
We have $H^0(S,\ell)=H^1(S,\ell)=0$ (Proposition \ref{q}), hence $\dim H^2(S,\ell)=\chi (\O_S)=2$. Thus the left hand side 
is $\Tr H^2(i_S,u)\in \{2,0,-2\}$. 
Since $f_0+f_1=40$ this gives $f_i\in \{16,20,24\}$; the case $f_0=24$ corresponds to $H^2(i_S,u)=1$.
 Moreover the number of fixed points of $i_X$ is $\#(B_2\cap X)=2f_0$. Thus the Proposition will follow if we can find $\ell$
in $\hat{A}_2$ with $f_0=24$.

\smallskip	
Put $F:=A_2\cap S$. Consider the homomorphism $\hat{A}_2\rightarrow (\Z/2)^F$ given by $\ell \mapsto (\langle a,\ell  \rangle)_{a\in F}$. For $ \ell \neq 0 $, the weight of the element $j(\ell)$ of $(\Z/2)^F$  (that is, the number of its nonzero coordinates) is   $f_1$, which 
 belongs to $\{16,20,24\}$. Therefore $ j $ is injective; its
 image is a 8-dimensional vector subspace of $(\Z/2)^F$, that is, a linear code, such that the weight of any nonzero vector belongs to $\{16,20,24\}$.
  A simple linear algebra lemma (\cite{Angers}, lemma 1) shows that a code in $(\Z/2)^{40}$ of dimension $\geq 7$ contains elements of weight $<20$; thus there exist elements $\ell$ in $\hat{A}_2$ with $f_1=16$, hence $f_0=24$.  \qed

\bigskip	

From now on we choose $\ell$ so that $i_S$ acts as trivially on $H^2(S,\ell)$. Thus  $i_X$  acts trivially on $H^0(X,K_X)$ and has $48$ fixed points. 
\begin{prop}
The canonical map $\fx : X\rightarrow \P^6$ factors through $i_X$, and induces an isomorphism of $X/i_X$ onto a complete intersection of $4$ quadrics in $\P^6$ with $48$ nodes.
\end{prop}
\pr Since $i_X$ acts as trivially on $H^0(X,K_X)$, we have a commutative diagram
\[\xymatrix @M=7pt{X\ar[d]^\pi \ar[r]^{\fx} &\Sigma \ar[d]^{p^{}_\Sigma } \ar@{^{(}->}[r]&\P^6 \ar@{-->}[d]^p\\
S\ar[r]^{\ft}& \Xi \ar@{^{(}->}[r]&\P^4
}\]where $\fx$ and $\ft$ are the canonical maps, $\Sigma $ and $\Xi$ their images,  $p$ the projection corresponding to the injection $H^0(S,K_S)\rightarrow H^0(X,K_X)$, $p_\Sigma $ its restriction to $\Sigma $.

The map $\ft\rond \pi :X\rightarrow \Xi$
 gives the quotient of $X$ by the action of $(\Z/2)^2$. Since $\tau $ acts non-trivially on $H^0(X,K_X)$, $\fx$ identifies $\Sigma $ with the quotient $X/i_X$. Thus all the maps in the left hand square of the above diagram are double coverings, \' etale outside finitely many points. In particular, since $K_X^2=32$, we have $\deg \Sigma =16$.

We choose bases  $(x_0,\ldots ,x_4)$  and $(u,v)$ of the $(+1)$ and $(-1)$-eigenspaces in $H^0(X,K_X)$ with respect to $\tau $.
The elements $u^2,uv,v^2$ of $H^0(X,K_X^{\otimes 2})$ are invariant under $\tau  $ and $i_X$, therefore they are pull-back of $i_S$-invariant forms in $H^0(S,K_S^{\otimes 2})$. Such a form comes from an element of $H^0(\Xi, \O_\Xi(2))$, hence from an element of $H^0(\P^4,\O_{\P}(2))$. 
Thus we have
 \[u^2=a(x)\qquad uv=b(x) \qquad v^2=c(x)\]
where $a,b,c$ are quadratic forms in $x_0,\ldots ,x_4$. Moreover the irreducible quadric $Q$ containing $\Xi$ is defined by a quadratic form $q(x)$ which vanishes on $\Sigma $. 

Thus $\Sigma $ is contained in  the subvariety $V$ of $\P^6$ defined by these 4 quadratic forms. If $V$ is a surface,  it has  degree $16$ and therefore is equal to $\Sigma $. Thus it suffices to prove that the morphism $p^{}_V: V\rightarrow Q$ induced by the projection $p$  is not surjective. 

Assume that $p^{}_V$ is surjective; it has degree 2, and we have a cartesian diagram
\[\xymatrix @M=7pt{\Sigma \ar[d]^{p^{}_\Sigma } \ar@{^{(}->}[r]& V \ar[d]^{p^{}_V}&\\
 \Xi \ar@{^{(}->}[r]&Q & \kern-2.3cm .
}\]

 The variety $V$ is irreducible: otherwise  $\Sigma $ is contained in one of its component, which maps birationally to $Q$, and  $p_{\Sigma }$  has degree 1, a contradiction.  
Since $Q\smallsetminus \mathrm{Sing}(Q)$ is simply connected, $p_V$ is branched along a surface $R\subset Q$. Since $\Xi$ is an ample divisor in $Q$ (cut out by a quartic equation), it meets $R$ along a curve, and $p_\Sigma $ is branched  along that curve, a contradiction.\qed

\begin{rem}
It follows that $\Xi=p(\Sigma )$ is defined by the equations $q(x)=b(x)^2-a(x)c(x)=0$. The 40 nodes of $\Xi$ break into two sets: the 16 points in $\P^4$ defined by $a(x)=b(x)=c(x)=q(x)=0$  are the images by $p_\Sigma $ of smooth points of $\Sigma $ fixed by the involution induced by $\tau $; 
 $p_\Sigma $ is \'etale over the other 24 nodes of $\Xi$, giving rise to the 48 nodes of $\Sigma $. 
\end{rem}

\begin{rem}
The two families of surfaces $ X $ that we have constructed are different; in fact, 
a surface $ X_1 $ of the first family is not even homeomorphic to a surface $ X_2 $ of the second one. Indeed $ X_1 $ admits an irrational genus 2 pencil $ X\rightarrow B $, and this is a topological property \cite{C}. But for a general member $ X_2 $ of the second family, 
the Albanese variety of the corresponding Schoen surface is simple \cite{S}, so its double cover $ \mathrm{Alb}(X_2)  $ is also simple; therefore $ X_2 $ cannot have an irrational  pencil of genus 2.

It follows that the corresponding surfaces $ \Sigma $ belong to two different connected components of the moduli space of complete intersections of 4 quadrics  in $\P^6$ with an even set of 48 nodes.

\end{rem}

\end{document}